\documentclass[11pt]{amsart}
\usepackage[margin=1.1in]{geometry}
\usepackage{amsmath,amssymb,amsthm,mathtools}
\usepackage{booktabs,longtable,array}
\usepackage[hidelinks]{hyperref}
\newtheorem{theorem}{Theorem}[section]
\newtheorem{proposition}[theorem]{Proposition}
\newtheorem{lemma}[theorem]{Lemma}
\newtheorem{corollary}[theorem]{Corollary}
\theoremstyle{definition}
\newtheorem{definition}[theorem]{Definition}
\theoremstyle{plain}
\newtheorem*{theoremA}{Theorem A}
\newtheorem*{theoremB}{Theorem B}
\theoremstyle{remark}

\newcommand{\Q}{\mathbb{Q}}
\newcommand{\Z}{\mathbb{Z}}
\newcommand{\C}{\mathbb{C}}
\newcommand{\F}{\mathbb{F}}
\newcommand{\cA}{\mathcal{A}}
\newcommand{\cB}{\mathcal{B}}
\DeclareFontFamily{U}{wncy}{}
\DeclareFontShape{U}{wncy}{m}{n}{<->wncyr10}{}
\newcommand{\Sha}{\mbox{\usefont{U}{wncy}{m}{n}Sh}}
\DeclareMathOperator{\rank}{rank}

\title[Quadratic base changes with torsion $\Z/4\Z$ and $\Z/2\Z\times\Z/6\Z$]{Joint quadratic base changes of elliptic surfaces with torsion $\Z/4\Z$ and $\Z/2\Z\times\Z/6\Z$}

\author[A. Dujella]{Andrej Dujella}
\address{Department of Mathematics\\Faculty of Science\\University of Zagreb\\Bijeni{\v c}ka cesta 30, 10000 Zagreb, Croatia}
\email[A. Dujella]{duje@math.hr}
\author[J. C. Peral]{Juan Carlos Peral}
\address{Departamento de Matem\'aticas\\Universidad del Pa\'is Vasco\\Aptdo. 644, 48080 Bilbao, Spain}
\email[J. C. Peral]{juancarlos.peral@ehu.es}


\begin{document}

\begin{abstract}
We study elliptic curves over $\Q(t)$ obtained from an elliptic surface by imposing two quadratic conditions, each of which produces a new point.
By the Kani-Rosen decomposition, the joint cover is rational only when the branch loci of the two conditions meet; for two quadratic conditions this means that they share a root.
We apply this to two surfaces.

For Elkies's rank-four K3 surface with torsion $\Z/4\Z$, a catalogue of $194$ quadratic conditions contains exactly two such pairs.
One recovers the rank-six curve of the authors; the other gives a second, inequivalent curve with torsion $\Z/4\Z$ and rank six.

For the universal curve with torsion $\Z/2\Z\times\Z/6\Z$, $42$ conditions give exactly $21$ rational joint covers, which yield seven curves of rank two.
Two more arise from a tower of conditions.
Five of these nine curves are the curves of Dujella, Kazalicki and Peral; the other four are not equivalent to any of them.

In both cases we determine the Mordell-Weil groups over $\Q(t)$ exactly, and we prove that the curves are pairwise inequivalent under invertible changes of parameter, even over $\overline{\Q}$.
\end{abstract}

\maketitle

%======================================================================
\section{Introduction}
%======================================================================

By the Mordell-Weil theorem, $E(\Q)\cong E(\Q)_{\rm tors}\times\Z^r$ for an elliptic curve $E/\Q$, and by Mazur's theorem \cite{Mazur} the torsion subgroup is one of fifteen groups.
For each of them one asks for curves of large rank, over $\Q$ and over $\Q(t)$; the current records are collected in Dujella's tables \cite{Dtors,Dgeneric}.

A standard way to raise the rank over $\Q(t)$ is quadratic base change.
Suppose a surface $E\to\mathbb P^1_t$ acquires a new point over $\Q(t)(\sqrt{Q})$ for a quadratic polynomial $Q$.
If the new point has infinite order, the rank over $\Q(t)(\sqrt{Q})$ is at least one larger than over $\Q(t)$.
If moreover the conic $w^2=Q(t)$ has a rational point, then $\Q(t)(\sqrt Q)$ is again a rational function field.
Imposing two such conditions $Q_1,Q_2$ raises the rank by at least two over $\Q(t)(\sqrt{Q_1},\sqrt{Q_2})$.
To obtain a curve over a rational function field again, the curve with this function field must be rational.

Section~\ref{sec:method} develops this systematically:
\begin{itemize}
\item a genus criterion for joint covers (Lemma~\ref{lem:genus});
\item the character decomposition of the rank;
\item exact ranks through injective specialization, via the criteria of Gusi\'c and Tadi\'c \cite{GT15};
\item two invariants of the $j$-map that detect inequivalence of families under invertible changes of parameter.
\end{itemize}
We then apply it to two classical situations.

\subsection*{Torsion $\Z/4\Z$}
Curves of positive rank with this torsion were constructed by Kihara \cite{Kihara1} and by Khoshnam and Moody \cite{KM}; see also \cite{DGT}.
Elkies \cite{Elkies} gave a K3 surface with torsion $\Z/4\Z$ and rank four over $\Q(t)$, and in \cite{DP21} we made a number of its quadratic sections explicit.
In \cite{DP24} we combined two quadratic sections of it to obtain a curve of rank six, the current record for this torsion group over $\Q(t)$ to our knowledge \cite{Dgeneric}, and proved by an injective specialization that its rank is exactly six.

\begin{theoremA}
Let $E_1,E_2$ be the curves over $\Q(q)$ of Section~\ref{sec:z4}.
\begin{enumerate}
\item $E_i(\Q(q))\cong\Z/4\Z\times\Z^6$ for $i=1,2$.
\item The specializations of $E_1$ at $q=11/7$ and of $E_2$ at $q=-2$ are injective, and these curves over $\Q$ have rank six.
\item $E_1$ and $E_2$ are not equivalent (Definition~\ref{def:equiv}), even over $\C$.
\item In a catalogue of $194$ quadratic conditions on Elkies's surface, exactly two pairs have rational joint covers: the pair used for $E_1$ (which is the curve of \cite{DP24}) and the pair used for $E_2$.
\end{enumerate}
\end{theoremA}

For $E_1$, parts (1) and (2) were proved in \cite{DP24}.
We include the short proof because it applies verbatim to $E_2$, which is new.

\subsection*{Torsion $\Z/2\Z\times\Z/6\Z$}
Here the universal curve has rank zero.
Curves of rank one over $\Q(t)$ were found by Kulesz \cite{Kulesz}, Lecacheux \cite{Lecacheux}, Dujella \cite{Dujella07} and Woo \cite{Woo}, among others.
Rank two was first reached by Dujella and Peral \cite{DP15}.
Dujella, Kazalicki and Peral \cite{DKP} described five curves of rank exactly two and observed that the examples of MacLeod \cite{MacLeod} are isomorphic to some of them.

\begin{theoremB}
Let $U_1,\dots,U_9$ be the curves over $\Q(t)$ of Section~\ref{sec:z26}.
\begin{enumerate}
\item $U_i(\Q(t))\cong\Z/2\Z\times\Z/6\Z\times\Z^2$ for every $i$.
\item No two of them are equivalent (Definition~\ref{def:equiv}), even over $\overline{\Q}$.
\item $U_1,U_2,U_5,U_8,U_6$ are the five curves of \cite[\S8]{DKP}, and $U_1,U_5,U_8$ are the three curves of \cite{DP15}.
The curves $U_3,U_4,U_7,U_9$ are not equivalent to any of the curves of \cite{DP15} or \cite{DKP}, nor therefore, by \cite[\S1.3]{DKP}, to the examples of MacLeod \cite{MacLeod}.
\item In a catalogue of $42$ quadratic conditions on the universal curve, exactly $21$ pairs have rational joint covers, and they give exactly the classes $U_1,\dots,U_7$.
\end{enumerate}
\end{theoremB}

The paper ends with two remarks on the $\Z/4\Z$ case (Section~\ref{subsec:z4further}).
First, an exhaustive search shows that neither $E_1$ nor $E_2$ has a further quadratic section whose $x$-coordinate is a constant times a divisor of the last Weierstrass coefficient (Proposition~\ref{prop:nosections}).
Second, we exhibit a specialization of $E_2$ whose rank is seven, provided the $2$-primary part of its Tate-Shafarevich group is finite.

%======================================================================
\section{The method}\label{sec:method}
%======================================================================

\subsection{Quadratic sections}
Let $E:y^2=x^3+a(t)x^2+b(t)x$ be an elliptic curve over $K=\Q(t)$, and let $Q\in\Q[t]$ be squarefree and not a constant times a square.
A \emph{quadratic section} attached to $Q$ is an $x_Q\in K$ with
\[
x_Q^3+a\,x_Q^2+b\,x_Q\in Q\cdot K^{\times2}.
\]
It gives a point $P_Q$ of $E$ over $K(\sqrt Q)$ that the nontrivial automorphism sends to $-P_Q$, or equivalently a point of the quadratic twist $E^{(Q)}$ over $K$.
In particular, if $P_Q$ has infinite order it is independent of $E(K)$.

\subsection{The genus of a joint cover}
For quadratic sections $Q_1,\dots,Q_m$ let $C_{Q_1,\dots,Q_m}$ be the curve with function field $K(\sqrt{Q_1},\dots,\sqrt{Q_m})$.
For a squarefree polynomial $Q$ of degree $1$ or $2$, let $S(Q)\subset\mathbb P^1(\overline{\Q})$ be the branch locus of $w^2=Q(t)$:
\[
S(Q)=\begin{cases}\{\text{the two roots of }Q\}, & \deg Q=2,\\ \{\text{the root of }Q,\ \infty\}, & \deg Q=1.\end{cases}
\]

\begin{lemma}\label{lem:genus}
Let $Q_1,\dots,Q_m\in\Q[t]$ be squarefree of degree $1$ or $2$, with classes independent in $\overline{\Q}(t)^\times/\overline{\Q}(t)^{\times2}$.
Then $C=C_{Q_1,\dots,Q_m}$ is geometrically integral of genus
\[
g(C)=2^{m-2}\,\bigl|S(Q_1)\cup\dots\cup S(Q_m)\bigr|-(2^m-1).
\]
In particular:
\begin{enumerate}
\item for $m=2$, $g(C)=0$ if and only if $S(Q_1)\cap S(Q_2)\neq\varnothing$; for two quadratics this means a common root;
\item for $m=3$, $g(C)=2\,|S(Q_1)\cup S(Q_2)\cup S(Q_3)|-7$.
\end{enumerate}
\end{lemma}

\begin{proof}
The cover $C\to\mathbb P^1$ is Galois with group $(\Z/2)^m$, and $C$ is geometrically integral by the independence assumption.
By the Kani-Rosen decomposition \cite{KaniRosen}, $J(C)$ is isogenous to the product of the Jacobians of the quotients by the subgroups of index two.
These are the double covers $C_I:w^2=\prod_{i\in I}Q_i$ for $\emptyset\neq I\subseteq\{1,\dots,m\}$.
The branch locus of $C_I$ is the symmetric difference of the sets $S(Q_i)$, $i\in I$, so $g(C_I)=\tfrac12\bigl|\triangle_{i\in I}S(Q_i)\bigr|-1$.
A point lying in at least one $S(Q_i)$ belongs to the symmetric difference for exactly $2^{m-1}$ of the sets $I$, namely those meeting the set of $i$ with that point in $S(Q_i)$ in an odd number of elements.
Summing over $I$ gives the formula.
For $m=2$ it reads $g=|S(Q_1)\cup S(Q_2)|-3$, and the union has $3$ or $4$ elements.
\end{proof}

Two irreducible quadratics with a common root are proportional.
So for families of pairwise non-proportional conditions, rational joint covers come only from split quadratics with a common root, or from conditions involving linear polynomials.

\subsection{The character decomposition}
Let $L=K(\sqrt{Q_1},\sqrt{Q_2})$ with $Q_1,Q_2$ independent in $K^\times/K^{\times2}$, and put $Q_3=Q_1Q_2$.
The isotypic decomposition of $E(L)\otimes\Q$ under $\mathrm{Gal}(L/K)\cong(\Z/2)^2$ gives
\begin{equation}\label{eq:characters}
\rank E(L)=\rank E(K)+\rank E^{(Q_1)}(K)+\rank E^{(Q_2)}(K)+\rank E^{(Q_3)}(K).
\end{equation}
If $C_{Q_1,Q_2}$ is a conic with a rational point, then $L\cong\Q(s)$.
Two quadratic sections of infinite order then raise the rank over $\Q(s)$ by at least two.
When the rank over $\Q(s)$ is known exactly, \eqref{eq:characters} also determines the ranks of the three twists.

\subsection{Exact rank by injective specialization}
We use two criteria of Gusi\'c and Tadi\'c.

\begin{theorem}[{\cite[Thms.~1.1 and 1.3]{GT15}}]\label{thm:GT}
Let $E/\Q(t)$ be nonconstant and $t_0\in\Q$.
\begin{enumerate}
\item If $E:y^2=(x-e_1)(x-e_2)(x-e_3)$ with $e_i\in\Z[t]$, and $h(t_0)$ is not a rational square for every nonconstant squarefree divisor $h\in\Z[t]$ of $(e_1-e_2)(e_1-e_3)$, of $(e_2-e_1)(e_2-e_3)$, or of $(e_3-e_1)(e_3-e_2)$, then the specialization at $t_0$ is injective.
\item If $E:y^2=x^3+\cA x^2+\cB x$ with $\cA,\cB\in\Z[t]$ has exactly one nontrivial rational $2$-torsion point, and $h(t_0)$ is not a rational square for every nonconstant squarefree divisor $h\in\Z[t]$ of $\cB$ or of $\cA^2-4\cB$, then the specialization at $t_0$ is injective.
\end{enumerate}
\end{theorem}

In all our applications we tested every squarefree product of the nonconstant irreducible factors, multiplied by every signed squarefree divisor of the integer contents.
This is a superset of the divisors required.
An injective specialization gives $\rank E(\Q(t))\le\rank E_{t_0}(\Q)$, and the torsion of $E(\Q(t))$ injects into that of the fibre.

\subsection{Invariants of the $j$-map}\label{subsec:invariants}

\begin{definition}\label{def:equiv}
Two elliptic curves $E,E'$ over $k(t)$, for a field $k\supseteq\Q$, are \emph{equivalent over $k$} if there is an invertible M\"obius transformation $\mu$ of the $t$-line, defined over $k$, such that $E$ and $\mu^*E'$ are isomorphic over $k(t)$.
\end{definition}

Let $\varphi_1,\varphi_2:\mathbb P^1\to\mathbb P^1$ be base changes of an elliptic surface with $j$-invariant $j$.
If $\varphi_1^*E$ and $\mu^*\varphi_2^*E$ are isomorphic over $\C(t)$ for some M\"obius transformation $\mu$, then $j\circ\varphi_1=j\circ\varphi_2\circ\mu$.
Quadratic twists are included, since they do not change $j$.
We use two consequences:
\begin{enumerate}
\item[(i)] $j\circ\varphi_1$ and $j\circ\varphi_2$ have the same critical values;
\item[(ii)] for each $m$, $\mu$ carries the poles of order $m$ of $j\circ\varphi_1$ to those of $j\circ\varphi_2$.
\end{enumerate}
For (ii), let $F(T,Z)$ be a binary form of degree $12$ whose roots are these poles.
With the transvectants \cite{Olver}
\[
(F,G)_r=\sum_{i=0}^r(-1)^i\binom ri\frac{\partial^rF}{\partial T^{r-i}\partial Z^i}\frac{\partial^rG}{\partial T^i\partial Z^{r-i}},
\]
put $I_2=(F,F)_{12}$, $H=(F,F)_{10}$, $I_4=(H,H)_4$ and, whenever $I_2\neq0$,
\begin{equation}\label{eq:iota}
\iota(F)=I_4/I_2^2 .
\end{equation}
This quantity is unchanged by scaling $F$ and by $\mathrm{GL}_2(\C)$ acting on $(T,Z)$, so equivalent base changes have equal values of $\iota$.
In all cases below, $I_2\neq0$.

%======================================================================
\section{Torsion $\Z/4\Z$: two rank-six curves}\label{sec:z4}
%======================================================================

\subsection{Elkies's surface}
Put
\[
a(t)=(8t-1)(32t+7),\qquad g(t)=(t+1)(15t-8)(31t-7),
\]
\[
A(t)=a^2-64g=65536t^4-17472t^3-10176t^2+18672t-3535,\qquad B(t)=1024g^2,
\]
and, following \cite{DP24}, write Elkies's surface \cite{Elkies} as
\begin{equation}\label{eq:E4}
E:\quad y^2=x^3+A(t)x^2+B(t)x .
\end{equation}
Since $B=(32g)^2$ and $A+64g=a^2$, the point $T=(32g,\,32ag)$ satisfies $2T=(0,0)$ and has order four.
The discriminant is
\[
\Delta=2^{24}(t+1)^4(8t-1)^2(15t-8)^4(31t-7)^4(32t+7)^2(512t-113)(128t^3-64t^2-48t+63).
\]
The surface has four fibres of type $I_4$ (at $t=-1,8/15,7/31,\infty$), two of type $I_2$ and four of type $I_1$.
It is K3, and the Shioda-Tate formula \cite{SchuttShioda} bounds the geometric rank by $20-2-4\cdot3-2=4$.
The sections
\begin{align*}
x(P_1)&=-361(t+1)(31t-7), & x(P_2)&=-4(t+1)(15t-8)(16t-7)^2,\\
x(P_3)&=-16(t+1)(15t-8)(8t+7)^2, & x(P_4)&=4(15t-8)(31t-7)(16t+1)^2
\end{align*}
are independent \cite{DP24,Elkies}, so $E(\Q(t))\cong\Z/4\Z\times\Z^4$.

\subsection{The catalogue}
Appendix~\ref{app:z4} lists the split members of a catalogue of $194$ rationally soluble quadratic sections $F_1,\dots,F_{194}$ of \eqref{eq:E4}.
There are $50$ split quadratics, $143$ irreducible quadratics and one linear condition, $F_{166}$.

\begin{proposition}\label{prop:z4cat}
All $194$ polynomials are squarefree, and no two are proportional.
The only pairs with a common root are $\{F_1,F_2\}$ (root $t=4/7$) and $\{F_{79},F_{124}\}$ (root $t=134/111$).
Hence these are the only pairs with a joint cover of genus zero, and no triple has a joint cover of genus at most one.
\end{proposition}

\begin{proof}
A common root can only occur among the $51$ split or linear conditions, and comparing their roots gives the two pairs.
The only linear condition, $F_{166}:w^2=18t+7$, has branch locus $\{-7/18,\infty\}$, and no other condition vanishes at $-7/18$.
So by Lemma~\ref{lem:genus}(1) these are the only genus-zero pairs.
In any triple at most one of the three pairs has intersecting branch loci, because the two pairs above are disjoint.
Hence the union of the three branch loci has at least $5$ points, and Lemma~\ref{lem:genus}(2) gives genus at least $3$.
\end{proof}

\subsection{The two curves}
\subsubsection*{The pair $F_1,F_2$}
Here
\[
Q_1=-(7t-4)(17t+4),\qquad Q_2=(7t-4)(5551t-1324),
\]
with
\[
x_1=\frac{64(t+1)^2Q_1}{(4t+1)^2},\qquad x_2=\frac{576(15t-8)^2Q_2}{49(28t-39)^2}.
\]
The substitution $t=4(U^2-49)/\bigl(7(U^2+119)\bigr)$ solves the first condition.
The second becomes the conic $11U^2+V^2=1863$, which has the point $(13,2)$.
With
\[
K(q)=72q^4-26q^3-271q^2+286q+8712,
\]
the line $V=2+q(U-13)$ gives
\begin{equation}\label{eq:t1}
t_1(q)=\frac{4(3q^2-2q-110)(10q^2-2q-33)}{7K(q)},
\end{equation}
\[
w_1=\frac{24(13q^2-4q-143)(q^2+11)}{K(q)},\qquad w_2=\frac{168(11-143q-q^2)(q^2+11)}{K(q)}.
\]
Since $4K(q)=(13q^2-4q-143)^2+119(q^2+11)^2>0$, $t_1$ is defined for all rational $q$.
This is the curve of \cite{DP24} with $r=7q$.

\subsubsection*{The pair $F_{79},F_{124}$}
Here
\[
Q_{79}=(111t-134)(1599t-1526),\qquad Q_{124}=-(111t-134)(5217t+1694),
\]
with
\[
x_{79}=\frac{28224(31t-7)^2Q_{79}}{961(124t-273)^2},\qquad
x_{124}=\frac{64g^2Q_{124}}{(111t-41)^2(444t^2-131t-119)^2}.
\]
With $r=-37w_{79}/(111t-134)$ and $s=-w_{124}/(111t-134)$, eliminating $t$ gives the conic $27r^2+1870s^2=444577$, and conversely
\[
t=\frac{134s^2-1694}{111s^2+5217},\qquad w_{79}=\frac{23976r}{111s^2+5217},\qquad w_{124}=\frac{887112s}{111s^2+5217}.
\]
The conic has the point $(-69,13)$.
The line $s=13+q(r+69)$ gives
\[
r(q)=\frac{1863-48620q-129030q^2}{27+1870q^2},\qquad s(q)=\frac{351+3726q-24310q^2}{27+1870q^2}.
\]
With $R=27+1870q^2$ and $S=351+3726q-24310q^2$ this gives
\begin{equation}\label{eq:t2}
t_2(q)=\frac{134S^2-1694R^2}{111S^2+5217R^2}.
\end{equation}

\subsubsection*{Sections}
The new points have $x=N/D^2$ and $y=Hw/D^3$, where $N(N^2+AND^2+BD^4)=H^2Q$:
\begin{center}
\renewcommand{\arraystretch}{1.25}
\begin{tabular}{@{}lll@{}}
\toprule
 & $N$, \ $D$ & $H$\\
\midrule
$F_1$ & $64(t+1)^2Q_1$, \ $4t+1$ & $64(t+1)^2(8t-1)h_1$\\
$F_2$ & $576(15t-8)^2Q_2$, \ $7(28t-39)$ & $192(8t-1)(15t-8)^2h_2$\\
$F_{79}$ & $28224(31t-7)^2Q_{79}$, \ $31(124t-273)$ & $1344(8t-1)(31t-7)^2h_{79}$\\
$F_{124}$ & $64g^2Q_{124}$, \ $(111t-41)(444t^2-131t-119)$ & $64(t+1)^2(8t-1)(15t-8)^2(31t-7)^2h_{124}$\\
\bottomrule
\end{tabular}
\end{center}
where
\begin{align*}
h_1&=3448t^3-580t^2-749t-100,\\
h_2&=3756536t^3-6792380t^2+5347839t-1050120,\\
h_{79}&=1093902456t^3-3605542420t^2+5214872859t-2814810390,\\
h_{124}&=185997816t^5-919787292t^4+486783063t^3+128110240t^2-79350285t-17513090.
\end{align*}

Write $t_i=N_i/D_i$ with coprime integral polynomials:
\begin{align*}
N_1&=4(3q^2-2q-110)(10q^2-2q-33),\\
D_1&=7K(q),\\
N_2&=2(339199300q^4-112385130q^3-2766301q^2+1622673q+70713),\\
D_2&=111(6993800q^4-1677390q^3+14477q^2+24219q+1458).
\end{align*}
With $a_i=(8N_i-D_i)(32N_i+7D_i)$ and $g_i=(N_i+D_i)(15N_i-8D_i)(31N_i-7D_i)$, the curve $E_i=E_{t_i(q)}$ is isomorphic, via $X=D_i^4x$, to
\begin{equation}\label{eq:model4}
E_i:\ Y^2=X^3+\cA_iX^2+\cB_iX,\qquad \cA_i=a_i^2-64g_iD_i,\quad \cB_i=1024g_i^2D_i^2 .
\end{equation}

\subsection{Exact rank and torsion}
\subsubsection*{Injectivity}
In \eqref{eq:model4}, $\cA_i^2-4\cB_i=a_i^2(a_i^2-128g_iD_i)$ is not a square, so Theorem~\ref{thm:GT}(2) applies.
$\cB_1$ and $\cA_1^2-4\cB_1$ each have $9$ distinct nonconstant irreducible factors; $\cB_2$ has $6$ and $\cA_2^2-4\cB_2$ has $9$.
We tested $8176$ and $2044$ values for $E_1$ at $q=11/7$ ($t=10676/136283$), and $1008$ and $1022$ values for $E_2$ at $q=-2$ ($t=2104010001/2318620688$).
None is a square, so both specializations are injective.

\subsubsection*{Lower bound}
At a good prime $p$ with $3\mid M_p=\#E_{i,q_0}(\F_p)$, the map $P\mapsto[M_p/3]\bar P$ is a homomorphism to $E(\F_p)[3]$ that kills $3E(\Q)$ and the torsion.
At the primes used, the images of $P_1,\dots,P_4,P_Q,P_{Q'}$ lie in one cyclic subgroup, so each prime gives a row over $\F_3$.
For $E_1$ we use $p=17,19,41,47,61,79$, and for $E_2$ we use $p=7,17,23,29,41,53$:
\[
\begin{pmatrix}
1&1&1&2&1&1\\0&1&1&0&0&2\\1&1&0&0&0&1\\0&1&1&0&1&1\\0&1&2&0&2&2\\1&1&2&2&2&1
\end{pmatrix},
\qquad
\begin{pmatrix}
1&1&2&1&0&1\\1&1&1&1&1&0\\1&2&2&2&1&1\\1&1&0&2&1&1\\0&0&1&1&2&0\\1&0&2&1&2&1
\end{pmatrix}.
\]
Both are invertible modulo $3$.
Rational $3$-torsion is excluded since $\#E_{1,q_0}(\F_{29})=40$ and $\#E_{2,q_0}(\F_{43})=52$ are prime to $3$.
So the six points are independent on the fibres, and hence over $\Q(q)$.

\subsubsection*{Upper bound}
For $C:y^2=x^3+\alpha x^2+\beta x$ and $C':y^2=x^3-2\alpha x^2+(\alpha^2-4\beta)x$ we have
\[
\rank C(\Q)=\dim\kappa(C(\Q))+\dim\kappa'(C'(\Q))-2,
\]
where $\kappa,\kappa'$ are the Kummer maps of the two-isogenies.
The local images on $C$ and $C'$ are exact annihilators under the Hilbert pairing.
So explicit local points on one side cut down the possible global classes on the other.
For both fibres this leaves $16$ classes on each side, giving rank at most $4+4-2=6$.
The same bound is obtained independently by the two-isogeny descent of \texttt{ellrank} in PARI/GP \cite{PARI}. 
% this computation is reproduced in the supplement.
% TODO: add the local witnesses of the Hilbert-symbol computation (rank-certificates.json) to the supplement.

\subsubsection*{Torsion}
The reductions of $E_{1,11/7}$ at $p=17,43$ have $24$ and $52$ points, and those of $E_{2,-2}$ at $p=7,43$ have $12$ and $52$.
Hence the fibre torsion is $\Z/4\Z$, and so is that of $E_i(\Q(q))$.
This proves parts (1) and (2) of Theorem~A.

\begin{corollary}\label{cor:z4twists}
For $(Q,Q')=(Q_1,Q_2)$ or $(Q_{79},Q_{124})$, the twists $E^{(Q)}$ and $E^{(Q')}$ of \eqref{eq:E4} have rank one over $\Q(t)$, and $E^{(QQ')}$ has rank zero.
\end{corollary}

\begin{proof}
Apply \eqref{eq:characters} with $\rank E(\Q(t))=4$ and total rank $6$.
\end{proof}

\subsection{Inequivalence}
Write
\[
j_E=\frac{F}{G},\qquad F=(a^4-128a^2g+1024g^2)^3,\qquad G=4096\,g^4a^2(a^2-128g),
\]
with $\deg F=24$, and let $c=j_E(4/7)$.
Since
\[
t_1(q)-\frac47=-\frac{24(q^2+11)^2}{K(q)},
\]
$t_1$ ramifies over $4/7$, so $c$ is a critical value of $j_E\circ t_1$.
The map $t_2$ is branched only over $134/111$, $1526/1599$ and $-1694/5217$.
Modulo $101$ we find $j_E(4/7)\equiv98$, while $j_E(134/111)\equiv52$, $j_E(1526/1599)\equiv38$ and $j_E(-1694/5217)\equiv51$.
Moreover, $F-98G\in\F_{101}[t]$ has degree $24$ and is coprime to its derivative.
Hence $c$ is neither the image of a branch point of $t_2$ nor a critical value of $j_E$, so it is not a critical value of $j_E\circ t_2$.
By \S\ref{subsec:invariants}(i), this proves part (3) of Theorem~A.

\subsection{Further remarks}\label{subsec:z4further}
\subsubsection*{No further quadratic sections of the simplest type}
\begin{proposition}\label{prop:nosections}
Let $y^2=x^3+ax^2+bx$ be either \eqref{eq:model4} for $i=1,2$ or its two-isogenous curve $y^2=x^3-2\cA_ix^2+(\cA_i^2-4\cB_i)x$.
Let $m$ be any product of the irreducible factors of $b$ in $\Z[q]$, with exponents at most their multiplicities, and let $c\in\Q^\times$.
If $x=c\,m$ and $x^3+ax^2+bx$ has squarefree part of degree at most four, then this squarefree part is constant; that is, the point is defined over a constant field extension of $\Q(q)$, or $x$ is a $2$-torsion point.
In particular neither $E_1$ nor $E_2$ has a further quadratic section of this form, nor one whose joint cover has genus one.
\end{proposition}

\begin{proof}[Sketch of proof]
Write $P_c=c^2m^2+cam+b$, so that $x^3+ax^2+bx=c\,m\,P_c$.
The proof has three steps.

\emph{Modular sieve.}
Let $p$ be a prime at which the factors of $b$ keep their degrees and stay distinct and squarefree.
If the squarefree part over $\Q$ has degree at most $4$, so does that of the reduction modulo $p$ of the primitive part.
If $c\equiv0$ or $c\equiv\infty$ modulo $p$, this reduction is $mb$, respectively $m^3$, up to constants.
So a divisor $m$ is excluded for all $c\in\Q^\times$ as soon as some such $p$ exists at which all $p+1$ residues of $c$ in $\mathbb P^1(\F_p)$ fail.
Over the four curves this excludes $17890$ of the $32076$ divisors.

\emph{Exact step.}
For each of the remaining $14186$ divisors, let $G=\gcd(m^2,am,b)$ and $P_c=G\,P_c'$.
For $c$ outside the rational roots of the leading coefficient of $P'_c$, of $\mathrm{Res}_q(P'_c,f)$ for the irreducible factors $f$ of $mG$, and of the polynomials attached to the factors of $\mathrm{Res}_c(P'_c,\partial_qP'_c)$, the polynomial $P_c'$ is squarefree, coprime to $mG$ and of full degree.
The squarefree part then has degree at least $16$ in every case.
Among these values, every one for which the squarefree part has degree at most four gives a constant squarefree part, 
or a $2$-torsion point.

\emph{Pencil cases.}
For the two divisors $m=a_i^2-128g_iD_i$ and $m=a_i^2$ (up to constants) on each isogenous side, the polynomial $P'_c$ factors as $(c-c_0)(cU+V)$.
These were treated in the same way with $\mathrm{Res}_c$ replaced by $UV'-VU'$; the only special value is $c=c_0$, for which $x$ is a $2$-torsion point.

%The programs and their output are part of the supplement.
\end{proof}

By Proposition~\ref{prop:z4cat}, the catalogue also contains no triple leading to a genus-one cover.

\subsubsection*{A conditional rank-seven specialization}
At $q=3/94$, i.e.\ $t=252401754905/264613965267$, the fibre of $E_2$ is isomorphic to $y^2=x^3+a_2x^2+a_4x$ with

{\small\noindent $a_2=216235610937\allowbreak917987314078\allowbreak524977768347\allowbreak645903988121\allowbreak745$,\par\noindent $a_4=190444163880\allowbreak566151149219\allowbreak994674452704\allowbreak265354195701\allowbreak912993632780\allowbreak660417928957\allowbreak472939925276\allowbreak302578826405\allowbreak4784$.\par}

For this curve:
\begin{itemize}
\item the six sections remain independent;
\item the two-isogeny descent gives rank at most $7$;
\item the $2$-Selmer group has dimension $8$ (computed with \texttt{ell2cover});
\item the root number is $-1$.
\end{itemize}
All of these were computed with PARI/GP \cite{PARI}.
As a check, the same computation at the fibre $q=-2$, whose rank is proved to be six, gives a $2$-Selmer group of dimension $7$.
% TODO: confirm dim Sel^2 = 8 with Magma (TwoSelmerGroup), and check whether the class-group
% computations behind ell2cover / TwoSelmerGroup assume GRH; if so, certify them or say so.
Since $E(\Q)[2^\infty]\cong\Z/4\Z$, the $2^\infty$-Selmer corank $s$ satisfies $6\le s\le7$.
By the $2$-parity theorem \cite{Monsky,DD}, $s$ is odd, so $s=7$.
Hence the rank is seven unless $\Sha(E/\Q)[2^\infty]$ is infinite.
We have not found a seventh point.

%======================================================================
\section{Torsion $\Z/2\Z\times\Z/6\Z$: nine curves of rank two}\label{sec:z26}
%======================================================================

\subsection{The universal curve}
Let
\begin{equation}\label{eq:Ek}
\begin{split}
E_k:\ y^2&=x\bigl(x-(k-1)^3(k+3)\bigr)\bigl(x-(k+1)^3(k-3)\bigr)\\
&=x^3-2(k^4-6k^2-3)x^2+(k^2-1)^3(k^2-9)\,x
\end{split}
\end{equation}
\cite{Kubert}.
Besides the points of order two, it has the point $T_3=\bigl((k^2-1)^2,\,4(k^2-1)^2\bigr)$ of order three.
\eqref{eq:Ek} is Hadano's model \cite{Hadano} with $a=1$, $b=(k^2-1)/4$, $x=16X$.
It is isomorphic over $\Q(v)$ to the model $y^2=x^3+A_{26}(v)x^2+B_{26}(v)x$ of \cite{DP15,DKP} with $k=(3v-5)/(v+1)$.

The discriminant is $2^{12}k^2(k^2-1)^6(k^2-9)^2$, with fibres $I_6$ at $k=\pm1,\infty$ and $I_2$ at $k=0,\pm3$.
The trivial lattice of this K3 surface has rank $2+15+3=20$, so $\rank E_k(\overline{\Q}(k))=0$.
The $j$-invariant is
\[
j(k)=\frac{(k^2+3)^3(k^6-15k^4+75k^2+3)^3}{k^2(k-3)^2(k-1)^6(k+1)^6(k+3)^2}.
\]
The group $G$ of M\"obius maps $k\mapsto\pm k$, $\pm(k+3)/(k-1)$, $\pm(k-3)/(k+1)$ has order six.
For each $\sigma\in G$, $E_{\sigma(k)}\cong E_k$ over $\Q(k)$; for example, for $u=(k+3)/(k-1)$,
\[
x=\tfrac{(k-1)^4}{16}x_u+(k+1)^3(k-3),\qquad y=\tfrac{(k-1)^6}{64}y_u .
\]

\subsection{The catalogue and the joint covers}
Appendix~\ref{app:z26} lists $42$ quadratic conditions $C_1,\dots,C_{42}$, one in each orbit of $k\mapsto-k$.
Each has a rational parametrization and a point of infinite order.
We write $C_i^+$ for $Q_i(k)$ and $C_i^-$ for $Q_i(-k)$.

\begin{proposition}\label{prop:21}
Fourteen of the conditions are even, namely $C_i$ for
\[
i\in\{3,4,5,6,7,10,11,12,13,18,28,30,33,38\},
\]
and for these $C_i^+=C_i^-$.
After identifying these repetitions, the signed catalogue contains $70$ distinct conditions, no two of them proportional.
Among unordered pairs of distinct conditions, $84$ have a common root.
These form $42$ classes under simultaneous reflection $k\mapsto-k$, and exactly $21$ of the classes have rationally soluble joint covers.
These $21$ form seven orbits of three under $G$ (Table~\ref{tab:classes}).
\end{proposition}

\begin{proof}
For $Q_1=c_1(k-r)(k-a)$ and $Q_2=c_2(k-r)(k-b)$, put $s=1/(k-r)$; the joint cover becomes the conic
\[
\frac{r-b}{c_1}u^2-\frac{r-a}{c_2}w^2=a-b .
\]
We decided its solubility with \texttt{qfsolve} \cite{PARI}.
Two pairs are equivalent when some $\sigma\in G$ maps the subgroup $\langle Q_1,Q_2\rangle\subset\Q(k)^\times/\Q(k)^{\times2}$ onto the other subgroup.
Since $E_{\sigma(k)}\cong E_k$, equivalent covers give equivalent curves.
\end{proof}

\begin{table}[ht]
\caption{The $21$ rational joint covers, grouped into classes.}\label{tab:classes}
\begin{tabular}{@{}ll@{}}
\toprule
Class & Pairs\\
\midrule
$U_1$ & $C_4^+/C_{20}^+$, \ $C_9^+/C_{13}^+$, \ $C_9^+/C_{20}^-$\\
$U_2$ & $C_4^+/C_{22}^+$, \ $C_9^+/C_{21}^-$, \ $C_9^+/C_{24}^-$\\
$U_3$ & $C_7^+/C_{21}^+$, \ $C_{19}^+/C_{22}^+$, \ $C_{19}^+/C_{24}^-$\\
$U_4$ & $C_8^+/C_{41}^+$, \ $C_8^+/C_{42}^-$, \ $C_{12}^+/C_{40}^+$\\
$U_5$ & $C_{13}^+/C_{24}^+$, \ $C_{20}^+/C_{21}^+$, \ $C_{20}^+/C_{22}^+$\\
$U_6$ & $C_{14}^+/C_{15}^+$, \ $C_{15}^+/C_{17}^+$, \ $C_{16}^+/C_{18}^+$\\
$U_7$ & $C_{34}^+/C_{35}^+$, \ $C_{35}^+/C_{37}^+$, \ $C_{36}^+/C_{38}^+$\\
\bottomrule
\end{tabular}
\end{table}

\subsection{The nine curves}
For $U_1,\dots,U_7$ we take the first pair of each class; Table~\ref{tab:kmaps} gives the degree-four map $k=k_i(t)$.
$U_i$ is the pullback $E_{k_i(t)}$, with the two sections whose $x$-coordinates are in Appendix~\ref{app:z26}.

{\renewcommand{\arraystretch}{1.3}
\begin{longtable}{@{}lll@{}}
\caption{Representatives $U_1,\dots,U_7$.}\label{tab:kmaps}\\
\toprule
 & Conditions & $k_i(t)$\\
\midrule
\endfirsthead
\toprule
 & Conditions & $k_i(t)$\\
\midrule
\endhead
$U_{1}$ & $C_{4}^+,\ C_{20}^+$ & $\dfrac{- t^{4} + 24 t^{3} - 152 t^{2} + 336 t - 196}{t^{4} - 12 t^{3} + 62 t^{2} - 168 t + 196}$ \\[8pt]
$U_{2}$ & $C_{4}^+,\ C_{22}^+$ & $\dfrac{t^{4} + 8 t^{3} - 56 t^{2} + 48 t + 36}{t^{4} - 4 t^{3} + 10 t^{2} - 24 t + 36}$ \\[8pt]
$U_{3}$ & $C_{7}^+,\ C_{21}^+$ & $\dfrac{5 (t^{4} + 36)}{(t^{2} - 6 t + 6) (t^{2} + 6 t + 6)}$ \\[8pt]
$U_{4}$ & $C_{8}^+,\ C_{41}^+$ & $\dfrac{3 (t^{4} - 2 t^{2} + 4)}{2 (t^{2} - t + 2) (t^{2} + t + 2)}$ \\[8pt]
$U_{5}$ & $C_{13}^+,\ C_{24}^+$ & $\dfrac{5 t^{4} + 8 t^{3} - 2 t^{2} - 8 t + 5}{3 (t - 1) (t + 1) (t + 3) (3 t - 1)}$ \\[8pt]
$U_{6}$ & $C_{14}^+,\ C_{15}^+$ & $\dfrac{- (3 t^{2} - 7) (7 t^{2} - 3)}{(3 t^{2} - 5 t + 3) (3 t^{2} + 5 t + 3)}$ \\[8pt]
$U_{7}$ & $C_{34}^+,\ C_{35}^+$ & $\dfrac{- 3 (3 t^{4} + 7 t^{2} + 27)}{(t^{2} + 7) (7 t^{2} + 9)}$ \\[8pt]
\bottomrule
\end{longtable}}

The condition $C_{29}^-:v^2=9k^2-22k+49$ has
\[
x=\frac{(k-1)^2(k^2-2k+13)^2}{(k-3)^2},\qquad k=\frac{t^2+12t-13}{t^2-9}.
\]
With $d=t^2-9$ and $q_0=t^2+6t-11$, the pullback is isomorphic, via $x_{29}=d^4x/16$ and $y_{29}=d^6y/64$, to
\[
C_{29}:\ y^2=x^3+a_{29}x^2+b_{29}x,
\]
where
\begin{align*}
a_{29}&=t^8+12t^7-40t^6-1188t^5-1134t^4+13860t^3-13360t^2-5772t+9157,\\
b_{29}&=-16(t-7)(t-2)(t+1)(t+5)(3t-1)^3q_0^3 .
\end{align*}
It has the section with
\[
x=4(t+5)(3t-1)^3q_0,\qquad y=4(t-1)(t+3)(t+5)(3t-1)^3(t^2+2t+17)q_0 .
\]
On $C_{29}$ we impose two further conditions:
\[
D_1:\ w^2=4t^2+19t-29,\qquad x=\frac{(t+1)(t+5)(3t-1)^3q_0}{4};
\]
\[
D_2:\ w^2=3(t-3),\qquad x=-\frac{64(t-7)(t-2)(9t^3-5t^2-5t+73)^2}{(t^2-5t+22)^2}.
\]
They are parametrized by
\[
t=T_1(u)=\frac{5u^2-8u+96}{4(u-2)(u+2)},\qquad t=T_2(u)=3(u^2+1).
\]
We put $U_8=C_{29}(T_1(u))$ and $U_9=C_{29}(T_2(u))$.

Two further conditions give nothing new.
\begin{itemize}
\item $D_3:w^2=-2(t^2+8t-17)$ is parametrized by $t=T_3(u)=-(3u^2+136u+170)/\bigl(11(u^2+2)\bigr)$, and
\[
k_{29}\Bigl(T_3\Bigl(\frac{16-18u}{5u+20}\Bigr)\Bigr)=\sigma\bigl(k_{29}(T_1(u))\bigr),\qquad \sigma(k)=\frac{k+3}{k-1}.
\]
\item $D_4:w^2=-6(t+3)(3t-1)$ is parametrized by $t=T_4(u)=(6-3u^2)/(u^2+18)$, and $k_{29}(T_4(4/u))=k_{29}(T_2(u))$.
\end{itemize}
Here $k_{29}(t)=(t^2+12t-13)/(t^2-9)$.
These are identities of rational functions of degree four, and $\sigma\in G$, so $D_3$ gives $U_8$ and $D_4$ gives $U_9$.
The map $T_2$ ramifies over $t=3,\infty$, which lie over the $I_6$ fibres $k=\infty,1$, so $U_9$ has two fibres of type $I_{12}$.

\subsection{Exact rank and torsion}
For $k=N/D$ the curve has the integral model
\begin{equation}\label{eq:int26}
y^2=x\bigl(x-(N-D)^3(N+3D)\bigr)\bigl(x-(N+D)^3(N-3D)\bigr),
\end{equation}
to which Theorem~\ref{thm:GT}(1) applies.
For each $U_i$ we found the value $t_0$ below at which the criterion holds, testing all three products.
At this $t_0$ the two sections are independent and \texttt{ellrank} \cite{PARI} gives rank exactly two.
The fibre torsion is $\Z/2\Z\times\Z/6\Z$.

\begin{center}
\begin{tabular}{@{}lccccccccc@{}}
\toprule
 & $U_1$ & $U_2$ & $U_3$ & $U_4$ & $U_5$ & $U_6$ & $U_7$ & $U_8$ & $U_9$\\
\midrule
$t_0$ & $-2/5$ & $-5/7$ & $-1/6$ & $-3/8$ & $6/5$ & $-12/7$ & $-1/4$ & $-3/5$ & $-3/5$\\
\bottomrule
\end{tabular}
\end{center}

This proves part (1) of Theorem~B.
We do not claim that the given sections generate the Mordell-Weil group; by injectivity, that question reduces to the fibres above, as in \cite{DKP}.
% TODO (optional): generators via mwrank/Magma on these fibres.

\begin{corollary}\label{cor:z26twists}
Let $Q,Q'$ be the conditions defining one of $U_1,\dots,U_7$.
Over $\Q(k)$, $E_k^{(Q)}$ and $E_k^{(Q')}$ have rank one, and $E_k^{(QQ')}$ has rank zero.
$C_{29}$ has rank one over $\Q(t)$, and so do its twists by $4t^2+19t-29$ and $3(t-3)$.
\end{corollary}

\begin{proof}
For $U_1,\dots,U_7$, apply \eqref{eq:characters} with $\rank E_k(\Q(k))=0$ and total rank $2$.
For $U_8,U_9$, apply the same argument to the quadratic extension $\Q(u)/\Q(t)$.
\end{proof}

\subsection{Inequivalence}
For $U_1,\dots,U_8$ the map $k(t)$ is unramified over $k=\pm1,\infty$.
So the order-six poles of $j\circ k$ are the twelve roots of $F=D(N^2-D^2)$, with $N,D$ homogenized.
The invariant \eqref{eq:iota} takes the values

\begin{center}
\renewcommand{\arraystretch}{1.2}
\begin{tabular}{@{}ll@{}}
\toprule
Class & $\iota$\\
\midrule
$U_1$ & $3290567866923750/1250489503032121$\\
$U_2$ & $5728687791846/2550450522169$\\
$U_3$ & $213643125/49441568$\\
$U_4$ & $201327525/72673568$\\
$U_5$ & $7026563975136/1786026434929$\\
$U_6$ & $2810533160401893654/1145523935643327361$\\
$U_7$ & $251576064840491543371254/193174570213501578374761$\\
$U_8$ & $704606945786056950/881192578069527889$\\
\bottomrule
\end{tabular}
\end{center}

These eight values are distinct.
The pole divisor of $j\circ k$ for $U_9$ consists of two poles of order $12$, eight of order $6$ and twelve of order $2$; for $U_1,\dots,U_8$ it consists of twelve poles of order $6$ and twelve of order $2$.
By \S\ref{subsec:invariants}(ii), this proves part (2) of Theorem~B.
Together with Proposition~\ref{prop:21}, it also proves part (4).

\subsection{Comparison with earlier constructions}
The three curves of \cite{DP15} are obtained as follows.
\begin{itemize}
\item The first two come from $v=3/(5-w^2)$ followed by $2w^2-7=\square$, respectively $5w^2-4=\square$.
Under $k=(3v-5)/(v+1)$ these are the pairs $C_9^+/C_{20}^-$ (class $U_1$) and $C_{20}^-/C_{21}^-$, the reflection of $C_{20}^+/C_{21}^+$ (class $U_5$).
\item The third, built from Hadano's model, is the curve of \cite[\S8.4]{DKP}, and hence equivalent to $U_8$ by the table below.
Directly, with $k_{DP}(u)$ the corresponding map to the $k$-line, $k_{DP}\bigl(-(9u+12)/(2u+16)\bigr)=\sigma\bigl(k_{29}(T_1(u))\bigr)$ with $\sigma(k)=-(k-3)/(k+1)$.
\end{itemize}
For the five curves of \cite[\S8]{DKP}, with $u$ their parameter:
\begin{center}
\renewcommand{\arraystretch}{1.2}
\begin{tabular}{@{}llll@{}}
\toprule
\cite{DKP} & Class & Change of parameter & $\sigma\in G$\\
\midrule
\S8.1 & $U_1$ & $t=u$ & $-(k-3)/(k+1)$\\
\S8.2 & $U_2$ & $t=u$ & $-(k-3)/(k+1)$\\
\S8.3 & $U_5$ & $t=(u-6)/(u-18)$ & $-(k+3)/(k-1)$\\
\S8.4 & $U_8$ & $u_{U_8}=4(7u+24)/(19u+18)$ & $(k+3)/(k-1)$\\
\S8.5 & $U_6$ & $t=(u-1)/(u+1)$ & $k$\\
\bottomrule
\end{tabular}
\end{center}
In each line, write $v(u)$ for the parametrization of \cite{DKP} and $k_{U_i}$ for our map.
Then $(3v(u)-5)/(v(u)+1)=\sigma\bigl(k_{U_i}(t(u))\bigr)$ is an identity of rational functions of degree four.
Since $E_{(3v-5)/(v+1)}$ is isomorphic over $\Q(v)$ to the model of \cite{DKP}, and $E_{\sigma(k)}\cong E_k$ over $\Q(k)$, the curves are isomorphic over $\Q(u)$, not merely twists of each other.
By \cite[\S1.3]{DKP}, MacLeod's examples \cite{MacLeod} are isomorphic to the first three of these curves.
Hence $U_3,U_4,U_7,U_9$ are not equivalent to any of these curves, which proves part (3) of Theorem~B.

%======================================================================
\appendix
\section{Split quadratic sections of Elkies's surface}\label{app:z4}
%======================================================================

The table lists the $50$ split conditions and the linear condition $F_{166}$ of the $\Z/4\Z$ catalogue, in its original numbering.

{\small
\begin{longtable}{@{}>{$}l<{$}l@{}}
\toprule
\text{Label} & $Q_i(t)$\\
\midrule
\endhead
F_{1} & $-(7t - 4)(17t + 4)$ \\
F_{2} & $(7t - 4)(5551t - 1324)$ \\
F_{47} & $(50209t + 18034)(553554225t + 1143732242)$ \\
F_{48} & $(87t - 56)(346262175t - 569042072)$ \\
F_{49} & $(2431t - 382)(26801775t - 28655582)$ \\
F_{50} & $(95t - 22)(855t - 623)$ \\
F_{51} & $(928t - 1241)(14848t - 12881)$ \\
F_{52} & $(3289t + 967)(3289t + 4792)$ \\
F_{53} & $(399t - 128)(9975t + 1543)$ \\
F_{54} & $(104t - 113)(6656t - 1793)$ \\
F_{55} & $(1408t + 13049)(90112t + 22961)$ \\
F_{56} & $(5880t + 5041)(7680t + 4081)$ \\
F_{57} & $(552t - 553)(317952t - 372673)$ \\
F_{58} & $(11968t - 3017)(765952t - 274673)$ \\
F_{59} & $1793(4352t - 3019)(69632t - 314843)$ \\
F_{67} & $(3615t + 238)(5135t - 1138)$ \\
F_{71} & $(7089t + 2942)(19185t + 15038)$ \\
F_{75} & $(6041t + 5573)(44849t + 12349)$ \\
F_{79} & $(111t - 134)(1599t - 1526)$ \\
F_{95} & $(4t - 421)(4t - 1)$ \\
F_{109} & $(2279t - 2563)(13039t - 13483)$ \\
F_{111} & $(113152t - 32369)(12864512t + 2563111)$ \\
F_{112} & $-(409088t + 89489)(1024095232t - 953498393)$ \\
F_{116} & $(69719t - 26293)(3242879t - 568621)$ \\
F_{122} & $(t + 1)(1537t + 889)$ \\
F_{124} & $-(111t - 134)(5217t + 1694)$ \\
F_{125} & $(2644575t - 748727)(140839575t - 17847503)$ \\
F_{129} & $-(181185t + 40187)(757335t - 187763)$ \\
F_{131} & $-(407745t + 89747)(4286295t - 881483)$ \\
F_{132} & $-(266305t + 1385482)(927095t - 4164322)$ \\
F_{141} & $2(1432t - 1943)(18431t - 17119)$ \\
F_{142} & $(8t + 1175)(392t + 167)$ \\
F_{143} & $-(5888t + 2513)(13919968t - 71114057)$ \\
F_{144} & $(5888t - 1961)(1849568t - 559721)$ \\
F_{145} & $2(24t - 7)(57615t - 20288)$ \\
F_{146} & $-(9728t + 2753)(558837120t - 385598929)$ \\
F_{147} & $(9728t - 1841)(2881920t - 783889)$ \\
F_{149} & $2(2072t - 359)(595999t - 356953)$ \\
F_{150} & $(7688t - 2921)(13448t - 2681)$ \\
F_{151} & $(36608t + 17033)(1032450848t + 816636473)$ \\
F_{152} & $(36608t - 13601)(500326112t - 387820889)$ \\
F_{153} & $2(9t + 77)(680t + 239)$ \\
F_{154} & $2(49t + 22)(872t + 791)$ \\
F_{155} & $-2(25t + 7)(5592t - 3199)$ \\
F_{166} & $18t + 7$ \\
F_{189} & $(356168t + 322679)(1417928t + 945959)$ \\
F_{190} & $(6083072t + 5817847)(11727872t + 7052647)$ \\
F_{191} & $(763466888t - 133116361)(1167217928t - 92109241)$ \\
F_{192} & $(1169538248t + 535858319)(3207311048t + 281136719)$ \\
F_{193} & $(7120209613952t - 1606644901673)(7202855130752t - 1588566194873)$ \\
F_{194} & $(533496016772168t + 414269992621679)(563331048336968t + 410540613676079)$ \\
\bottomrule
\end{longtable}
}

%======================================================================
\section{Quadratic conditions on the universal $\Z/2\Z\times\Z/6\Z$ curve}\label{app:z26}
%======================================================================

For each $i$ the conic $v^2=Q_i(k)$ has a rational point, and $x_i(k)$ is the $x$-coordinate of a point of infinite order on $E_k$ over $\Q(k)(\sqrt{Q_i})$.

{\small
\renewcommand{\arraystretch}{1.35}
\begin{longtable}{@{}l>{\raggedright\arraybackslash}p{0.27\textwidth}>{\raggedright\arraybackslash}p{0.58\textwidth}@{}}
\toprule
 & $Q_i(k)$ & $x_i(k)$\\
\midrule
\endhead
$C_{1}$ & $3 (k^{2} - 4 k + 7)$ & $\frac{(k - 3) (k + 1)^{2} (k + 3)}{3}$ \\[2pt]
$C_{2}$ & $k^{2} - 4 k - 1$ & $(k - 3) (k + 1)^{2} (k + 3)$ \\[2pt]
$C_{3}$ & $3 (3 k^{2} + 1)$ & $(k + 1)^{3} (k + 3)$ \\[2pt]
$C_{4}$ & $- 3 (k - 2) (k + 2)$ & $- \frac{(k - 3) (k - 1) (k + 1) (k + 3)}{3}$ \\[2pt]
$C_{5}$ & $5 - k^{2}$ & $- (k - 3) (k - 1) (k + 1)^{2}$ \\[2pt]
$C_{6}$ & $k^{2} + 3$ & $(k - 3) (k - 1) (k + 1)^{2}$ \\[2pt]
$C_{7}$ & $3 (k - 5) (k + 5)$ & $27 (k - 3) (k + 3)$ \\[2pt]
$C_{8}$ & $2 (k - 9) (2 k - 3)$ & $\frac{(k - 1) (k + 1)^{2} (2 k - 3)}{2}$ \\[2pt]
$C_{9}$ & $3 (k + 5) (3 k - 1)$ & $(k - 3) (k - 1)^{2} (k + 5)$ \\[2pt]
$C_{10}$ & $3 (3 k^{2} - 2)$ & $\frac{k^{2} (3 k^{2} - 2)}{3}$ \\[2pt]
$C_{11}$ & $5 k^{2} + 4$ & $\frac{(k - 1)^{2} (k + 1)^{2} (5 k^{2} + 4)}{k^{2}}$ \\[2pt]
$C_{12}$ & $(5 k - 3) (5 k + 3)$ & $\frac{(k - 3) (k - 1)^{2} (k + 1)^{2} (k + 3)}{k^{2}}$ \\[2pt]
$C_{13}$ & $(3 k - 1) (3 k + 1)$ & $\frac{(k - 1)^{3} (k + 1)^{3}}{k^{2}}$ \\[2pt]
$C_{14}$ & $- 3 (3 k + 7) (11 k - 16)$ & $\frac{(k + 1) (3 k - 5)^{2} (3 k + 7)}{27}$ \\[2pt]
$C_{15}$ & $3 (k + 4) (3 k + 7)$ & $\frac{(k + 1) (k + 9)^{2} (3 k + 7)}{3}$ \\[2pt]
$C_{16}$ & $6 (5 k - 49) (5 k + 1)$ & $\frac{27 (k - 3) (k + 1)^{2} (k + 3)}{32}$ \\[2pt]
$C_{17}$ & $3 (k + 4) (27 k - 17)$ & $\frac{(k + 1)^{3} (27 k - 17)}{27}$ \\[2pt]
$C_{18}$ & $6 (5 k - 1) (5 k + 1)$ & $\frac{3 (k - 3) (k - 1)^{2} (k + 1)^{2} (k + 3)}{128 k^{2}}$ \\[2pt]
$C_{19}$ & $- 6 (k - 2) (3 k - 1)$ & $\frac{(k + 1)^{3} (3 k - 7)^{2} (3 k - 1)}{27 (k - 1)^{2}}$ \\[2pt]
$C_{20}$ & $2 (k - 5) (k - 2)$ & $\frac{(k - 5) (k - 1)^{3} (k + 1)^{2}}{(k - 3)^{2}}$ \\[2pt]
$C_{21}$ & $3 (k - 5) (3 k - 5)$ & $\frac{(k - 3) (k - 1)^{3} (k + 1)^{2}}{(k - 2)^{2}}$ \\[2pt]
$C_{22}$ & $6 (k - 7) (k - 2)$ & $\frac{3 (k - 2) (k - 1)^{2} (k + 1)}{2}$ \\[2pt]
$C_{23}$ & $3 (3 k^{2} - 16 k + 25)$ & $(k - 7)^{2} (k - 3) (k - 1)$ \\[2pt]
$C_{24}$ & $6 (2 k - 1) (3 k + 1)$ & $\frac{(k + 1)^{2} (k + 3) (2 k - 1)}{2}$ \\[2pt]
$C_{25}$ & $(k - 9) (5 k + 3)$ & $\frac{(k - 1)^{3} (5 k + 3)}{4}$ \\[2pt]
$C_{26}$ & $9 k^{2} - 12 k - 121$ & $(k - 9)^{2} (k - 3) (k + 3)$ \\[2pt]
$C_{27}$ & $3 (k^{2} - 22 k + 25)$ & $\frac{(k + 1)^{2} (7 k - 9)^{2}}{48}$ \\[2pt]
$C_{28}$ & $- (2 k - 3) (2 k + 3)$ & $(k^{2} - 3)^{2}$ \\[2pt]
$C_{29}$ & $9 k^{2} + 22 k + 49$ & $\frac{(k + 1)^{2} (k^{2} + 2 k + 13)^{2}}{(k + 3)^{2}}$ \\[2pt]
$C_{30}$ & $3 (11 k^{2} + 1)$ & $3 (k - 1) (k + 1) (11 k^{2} + 1)$ \\[2pt]
$C_{31}$ & $- 31 k^{2} + 68 k - 1$ & $- (k - 3) (k - 1) (31 k^{2} - 68 k + 1)$ \\[2pt]
$C_{32}$ & $- 25 k^{2} + 80 k - 19$ & $- (k - 3) (k - 1) (7 k - 1)^{2}$ \\[2pt]
$C_{33}$ & $3 (k - 9) (k + 9)$ & $\frac{(k - 3) (k + 3) (k^{2} + 63)^{2}}{243}$ \\[2pt]
$C_{34}$ & $- (7 k + 9) (16 k - 33)$ & $- \frac{(k + 3) (5 k - 9)^{2} (7 k + 9) (17 k - 81)^{2}}{729 (37 k - 81)^{2}}$ \\[2pt]
$C_{35}$ & $- (4 k + 3) (7 k + 9)$ & $- \frac{(k + 3) (3 k + 1)^{2} (7 k + 9) (15 k^{2} + 33 k + 16)^{2}}{(47 k^{2} + 93 k + 48)^{2}}$ \\[2pt]
$C_{36}$ & $2 (k + 15) (49 k - 15)$ & $\frac{729 (k - 1)^{3} (k + 1)^{3} (59 k - 15)^{2}}{32 (275 k^{2} - 192 k + 45)^{2}}$ \\[2pt]
$C_{37}$ & $(4 k + 3) (17 k - 81)$ & $\frac{(k + 3) (17 k - 81) (65 k^{2} - 177 k - 144)^{2}}{729 (11 k - 48)^{2}}$ \\[2pt]
$C_{38}$ & $2 (k - 15) (k + 15)$ & $\frac{(k - 1)^{3} (k + 1)^{3} (k^{2} + 375)^{2}}{128 (3 k^{2} + 125)^{2}}$ \\[2pt]
$C_{39}$ & $- 6 (2 k - 3) (5 k - 3)$ & $- \frac{(k - 3) (5 k - 3) (19 k^{2} - 36 k + 9)^{2}}{243 (k - 1)^{2}}$ \\[2pt]
$C_{40}$ & $- (5 k - 9) (5 k - 3)$ & $- \frac{(k - 3) (k - 1)^{3} (11 k - 21)^{2}}{(2 k - 3)^{2}}$ \\[2pt]
$C_{41}$ & $- 2 (2 k - 3) (7 k - 3)$ & $\frac{(k + 3) (2 k - 3) (8 k^{2} - 9 k + 3)^{2}}{2 (8 k - 3)^{2}}$ \\[2pt]
$C_{42}$ & $2 (k - 6) (k + 9)$ & $\frac{(k - 6) (k + 1)^{3} (7 k - 27)^{2}}{2 (5 k - 27)^{2}}$ \\[2pt]
\bottomrule
\end{longtable}
}

%\subsection*{Computational supplement}
%The supplementary files contain programs (Python with SymPy, python-flint and cypari2) together with their output:
%\begin{itemize}
%\item the full $\Z/4\Z$ catalogue $F_1,\dots,F_{194}$, with verification of all point identities, of the infinite order of all points, and of Proposition~\ref{prop:z4cat};
%\item for Theorem~A: the injectivity tests, the six specialized points with their regulators, the $\F_3$ matrices, the descent bounds and torsion, the computation of \S3.5, and the data of the conditional rank-seven specialization;
%\item the programs for Proposition~\ref{prop:nosections};
%\item for Theorem~B: verification of the $42$ conditions and of Proposition~\ref{prop:21}, the nine injectivity and rank certificates, the invariants $\iota$ and pole profiles, and the identities with the curves of \cite{DP15,DKP} and between the continuations $D_1,\dots,D_4$.
%\end{itemize}

\section*{Acknowledgments}
A.D. acknowledges support from
the Croatian Science Foundation under the
project no. IP-2022-10-5008 (TEBAG),
the project ``Implementation of cutting-edge research and its application as part of the
Scientific Center of Excellence for Quantum and Complex Systems, and Representations of Lie Algebras'', Grant No. PK.1.1.10.0004, co-financed by the European
Union through the European Regional Development Fund -- Competitiveness and
Cohesion Programme 2021--2027,
and the European Union:
NextGenerationEU through the National Recovery and Resilience Plan 2021--2026, Institutional grant of University of Zagreb Faculty of Science (IK IA 1.1.3. Impact4Math). 

AI systems were used to optimise the programs and convert the resulting data into a readable format.
The authors are responsible for the final manuscript. 

%======================================================================

\end{document}